\documentclass[10pt,a4paper]{smfart}
\usepackage[all]{xy}
\usepackage{a4}
\usepackage{microtype}
\usepackage[french]{babel}

\usepackage[T1]{fontenc}

\usepackage{a4,amssymb,amsthm,amsmath}
\usepackage{lmodern}

\def\à{\`a}
\def\è{\`e}
\def\ä{\"a}

\newcommand{\R}{\mathbb{R}}
\newcommand{\C}{\mathbb{C}}

\newcommand{\X}{\mathcal{X}}
\newcommand{\Y}{\mathcal{Y}}

\newcommand{\T}{\mathbb{T}}

\newcommand{\V}{\mathcal{V}}
\newcommand{\NN}{\mathcal{N}}
\newcommand{\HH}{\mathcal{H}}

\newcommand{\Sp}{\mathbb{S}}

\newcommand{\ZZ}{\mathcal{Z}}
\newcommand{\DD}{\mathcal{D}}
\newcommand{\OO}{\mathcal{O}}
\newcommand{\RR}{\mathcal{R}}

\newcommand{\J}{\mathcal{J}}
\newcommand{\JJ}{\mathbb{J}}

\newcommand{\lra}{\longrightarrow}
\newcommand{\lms}{\longmapsto}
\newcommand{\bw}{\bigwedge}
\newcommand{\w}{\wedge}
\newcommand{\wh}{\widehat}

\theoremstyle{definition}

\begin{document}
\frontmatter
 \title{Espace de twisteurs des structures complexes généralisées}
\author{Guillaume~Deschamps}
\address{Université de Brest}
\email{Guillaume.Deschamps@univ-brest.fr}
\date{\today}

\begin{abstract} Le but de cet article est d'utiliser les
structures complexes généralisées pour étendre la définition
d'espace de twisteurs introduite par Penrose. Ainsi à toute
4-variété riemannienne $(M,g)$ nous associons  le fibré
$\ZZ(M,g)\lra M$ des structures presque complexes généralisées sur
$M$ compatibles avec $g$. Comme dans l'article d'Atiyah, Hitchin
et Singer \cite{AHS78}, nous verrons que $\ZZ(M,g)$ admet une
structure presque complexe généralisée $\JJ$ dont nous donnerons
un critère d'intégrabilité. Ceci permettra de  construire une
passerelle entre la géométrie riemannienne sur $(M,g)$ et la
géométrie complexe généralisée sur $\ZZ(M,g)$. Dans une dernière
partie nous verrons comment étendre ces constructions aux variétés
quaternioniques Kähler et ferons  le lien avec un résultat de
Bredthauer \cite{Bre}.
\end{abstract}

\maketitle \setcounter{tocdepth}{4} \tableofcontents

\mainmatter
\section{Introduction}
 Le concept de structure complexe généralisée a été introduit
par Hitchin \cite{Hit} dans le but d'unifier les notions de
structure presque complexe et de structure presque symplectique.
C'est ensuite Gualtieri \cite{Gua2} qui a donné son essor à cette
théorie qui intéresse désormais tout autant les mathématiciens que
les physiciens. On trouvera des références dans \cite{DM1}.

Introduite par Penrose \cite{Pen},  la théorie des twisteurs
permet de son côté d'associer à toute 4-variété riemannienne
$(M,g)$ un fibré $Z(M,g)\lra M$  en sphères  $\Sp^2$ dont l'espace
total est l'espace des structures presque complexes sur $M$
compatibles avec $g$.  Un des attraits de cette théorie est de
coder des propriétés géométriques de $(M,g)$ en termes de
structure holomorphe sur son espace de twisteurs. En particulier,
on montre que $Z(M,g)$ admet une structure presque complexe
canonique $\JJ$ dont l'intégrabilité  dépend de la courbure de $g$
\cite{AHS78}.

Le but de cet article est d'étendre la construction de Penrose aux
structures complexes généralisées, c'est-à-dire d'étudier le fibré
des twisteurs $\ZZ(M,g)\lra M$ défini comme  le fibré des
structures presque complexes généralisées sur $M$ compatibles avec
la métrique $g$. Nous verrons que l'espace total admet une
structure presque complexe généralisée $\mathbb J$ naturelle et
que les fibres admettent quatre composantes connexes isomorphes à
$\Sp^2\times\Sp^2$. Nous donnerons un critère d'intégrabilité de
$\mathbb J$ qui dépendra bien entendu de la composante connexe
considérée. Cela nous donnera une caractérisation simple des
métriques d'Einstein, des métriques Ricci plates anti-autoduales
ou des métriques à courbure sectionnelle constante en terme
d'intégrabilité d'une structure complexe généralisée sur
$\ZZ(M,g)$ (cf. théorèmes 1 \& 2).

En dimension strictement plus grande que 4, et pour toute
4$n$-variété riemannienne $(M,g)$ munie d'une structure
quaternionique Kähler $D$, Salamon \cite{Sal} a défini le fibré
des twisteurs $Z(M,D)\lra M$ comme le fibré des structures presque
complexes de $D$. C'est un fibré en sphères $\Sp^2$ muni d'une
structure presque complexe canonique $\JJ$ qui est automatiquement
intégrable. Dans la dernière partie de ce papier nous verrons
qu'on peut encore associer à toute variété munie de deux
structures quaternioniques Kähler $(M,g,D_1,D_2)$, un fibré des
twisteurs $\ZZ(M,D_1,D_2)\lra M$. C'est un fibré de fibres
$\Sp^2\times\Sp^2$ muni d'une structure presque complexe
généralisée $\JJ$ dont nous étudierons l'intégrabilité. Nous
verrons en quoi c'est la généralisation naturelle de la situation
en dimension quatre et en quoi cela étend un résultat de
Bredthauer établi sur les variétés hyperkählériennes généralisées
\cite{Bre}.

\quad\\

\section{Préliminaires}
\subsection{Structures presque complexes généralisées.}
 Soit $M$ une
variété  de dimension $2n$. En géométrie généralisée on étudie non
pas le fibré  tangent
  de $M$ noté $TM$ mais plutôt la somme du fibré tangent et du fibré cotangent
   que  nous noterons $\T M=TM\oplus T^\star M$. Sur $\T M$ on a
une pseudo-métrique naturelle de signature $(2n,2n)$ définie par :
$$
<X+\xi,Y+\eta>=\frac{1}{2}\Big(\xi(Y)+\eta(X)\Big), \quad\forall
X,Y\in TM\textrm{ et }\forall\xi,\eta\in  T^\star M.
$$
Une structure {\it presque complexe} sur $M$ est la donnée d'un
endomorphisme $J$ de $TM$ tel que $J^2=-Id$. Une structure {\it
presque symplectique}
 sur $M$ est la donnée d'une $2$-forme anti-symétrique
  non dégénérée $w\in \bigwedge^2T^\star M$. En utilisant le
produit intérieur, on peut voir $w$ comme une application
$w:TM\lra T^\star M$ telle que $w^\star=-w$, où $w^\star$ est
l'adjoint de $w$. La motivation première de la géométrie
généralisée est d'unifier ces deux notions.

\quad\\
{\bf Définition \cite{Hit,Gua2}.} Une structure {\it presque
complexe généralisée} sur $M$ est la donnée d'un endomorphisme
$\mathcal J$ sur $\T M$ qui vérifie d'une part que $\mathcal J$
est presque complexe: $\mathcal J^2=-Id$ et d'autre part que
$\mathcal J$ est presque symplectique : $\mathcal
J^\star=-\mathcal J$. Ou de manière équivalente, une structure
presque complexe généralisée sur $M$ est la donnée d'une structure
presque complexe sur $\T M$ orthogonale pour la pseudo-métrique
définie précédemment.

\quad\\
{\bf Remarque.} On peut montrer qu'il existe une structure presque
complexe
 généra\-lisée sur $M$ seulement si
 sa dimension est paire et si $M$ satisfait
  certaines conditions topologiques \cite{Gua2}.

\quad\\
{\bf Propriété \cite{Hit,Gua2}.} Une structure presque complexe
généralisée $\J$ sur $M$ équivaut à la donnée d'un champ de
sous-espaces isotropes maximaux $L\subset \T M\otimes \C$ tel que
$L\cap \overline L=\{0\}$.

\quad\\
 On note $pr_1:(TM\oplus T^\star M)\otimes \C\lra TM\otimes\C$ la
première projection.

\quad\\
{\bf Définition \cite{Hit,Gua2}.} La codimension de  $pr_1(L)$
dans $TM\otimes \C$ est un invariant de la structure presque
complexe généralisée appelé le {\it type} de $\J$.

\quad\\
Comme le montre les exemples suivants, la notion de structure
complexe généralisée regroupe sous un même formalisme les notions
de structure complexe et de structure symplectique.

\quad\\
{\bf Exemple 1.}  Une structure presque complexe $J$ sur $M$
définit la structure presque complexe généralisée $\mathcal
J_J=\left(\begin{array}{cc} J&0\\0&-J^\star\end{array}\right)$ où
$J^\star$ est l'adjoint de $J$ sur $T^\star M$.  Le type de $\J_J$
est $n$.

\quad\\
{\bf Exemple 2.} De même une structure presque symplectique $w$
sur $M$ définit la structure presque complexe généralisée
$\mathcal J_w=\left(\begin{array}{cc}
0&-w^{-1}\\w&0\end{array}\right)$. Le type de $\J_w$ est $0$.

\quad\\
{\bf Exemple 3.} Toute 2-forme différentielle $B$ sur $M$ définit
l'application orthogonale
$$
\begin{array}{lccc}
e^B :&TM\oplus T^\star M&\lra&TM\oplus T^\star M\\
&X+\xi&\lms&X+\xi+i_X B
\end{array}
$$
où $i_X$ est le {\it produit intérieur}. Si $\J$ est une structure
presque complexe généralisée sur $M$, alors $e^{-B}\J
e^B$ aussi. Une telle transformation préserve le type.\\

\subsection{Intégrabilité.}
 On note
$[X,Y]$ le {\it crochet de Lie} de deux champs de vecteurs $X, Y$
sur $M$
 et  $\mathcal L_X$ la {\it dérivée de Lie} suivant le champ $X$.
 On définit  le {\it crochet de Courant} \cite{Cou} pour tout
$X+\xi,Y+\eta\in \mathcal C^\infty(\T M)$ par
$$[X+\xi,Y+\eta]=[X,Y]+\mathcal L_X\eta-\mathcal L_Y\xi-\frac{1}{2}d(i_X\eta-i_Y\xi).$$

\quad\\
{\bf Remarque.} Les notations ne sont pas ambiguës  dans la mesure
où le crochet de Courant et le crochet de Lie coïncident  sur
 les champs de vecteurs. Par
contre le crochet de Courant ne vérifie pas l'identité de Jacobi.

\quad\\
{\bf Définition \cite{Hit,Gua2}.} Une structure presque complexe
généralisée $\J$ est dite {\it intégrable} si le tenseur de
Nijenhuis $\NN$ défini par :
$$
\NN(\mathcal X, \mathcal Y)=[\J \mathcal X,\J \mathcal Y]-\J[\J
\mathcal X,\mathcal Y]-\J[ \mathcal X,\J \mathcal Y]-[\mathcal
X,\mathcal Y]\qquad\forall \mathcal X,\mathcal Y\in C^\infty(\T M)
$$
est nul sur $\T M$. En terme de champ de sous-espaces isotropes
maximaux $L$ cela équivaut à demander à l'espace des sections de
$L$ d'être stable par crochet de Courant.

\quad\\
 Cette définition d'intégrabilité est naturelle
au sens où elle généralise la notion d'intégrabilité des
structures presque complexes et presque symplectiques comme le
montre la proposition suivante.

\quad\\
{\bf Propriété \cite{Gua2}.}
\begin{enumerate}
\item[a)] Une structure presque complexe $J$ sur $M$ est
intégrable si et seulement si la structure presque complexe
généralisée associée $\J_J$ est intégrable.

\item[b)]  Une structure presque symplectique $w$ sur $M$ est
intégrable (c'est-à-dire $w$ fermée)  si et seulement si $\J_w$
est intégrable.

\item[c)] Soient $\J$ une structure presque complexe généralisée
et $B$ une 2-forme différentielle fermée sur $M$. Alors $\J$  est
intégrable si et seulement si sa $B$ transformation $e^{-B}\J e^B$
l'est.

\item[d)] Une variété sans aucune structure complexe et sans
aucune structure symplectique peut admettre une structure complexe
généralisée \cite{GC}.
\end{enumerate}

\subsection{Espaces des twisteurs "classiques" en dimension 4}

 On considère ici $(M,g)$ une
$4$-variété riemannienne orientée connexe et on note $O_g(TM)$ le
fibré des endomorphismes de $TM$ orthogonaux pour la métrique $g$.

\quad\\
{\bf Définition \cite{AHS78}.} Le fibré des twisteurs d'une
4-variété riemannienne $(M,g)$ défini par Atiyah, Hitchin et
Singer est le fibré $\pi : Z(M,g)\lra M$ des structures presque
complexes sur $M$ compatibles avec la métrique $g$:
$$
Z(M,g)=\{u\in O_g(TM)/u^2=-Id\}.
$$
\quad\\
 La fibre $Z(\R^4):=\{u\in O(4)/u^2=-Id\}$ est difféomorphe
à $O(4)/U(2)$ et admet donc deux composantes connexes, chacune
isomorphe à $\Sp^2$. Notons que, si $u\in O(4)$ vérifie $u^2=-Id$,
alors nécessairement $u\in SO(4)$. On dira donc qu'un
endomorphisme $u$ de $ TM$ est {\it compatible avec
l'orientation}, et on notera $u>>0$ si pour tous champs de
vecteurs $(X, Y)$, la famille $(X,uX, Y, uY)$ est linéairement
dépendante ou positivement orientée. Les deux composantes connexes
de $Z(M,g)$ sont
$$
Z^+(M,g)=\{u\in O_g(TM)/u^2=-Id\textrm{ et }u>>0\} ,
$$
$$
Z^-(M,g)=\{u\in O_g(TM)/u^2=-Id\textrm{ et }u<<0\}.
$$
Ces deux fibrés en sphères ont pour groupe structural $SO(3)$.
Leurs fibres peuvent donc être munies de la structure complexe de
$\C P^1$.

 Plus généralement, sur $Z(M,g)$ on peut
définir une structure presque complexe naturelle. En effet, la
connexion de Levi-Civita induit une décomposition du fibré tangent
$TZ(M,g)=\HH\oplus \V$ en la somme d'une
 distribution horizontale $\HH$ et d'une distribution verticale  (i.e. tangente aux fibres
$\V=ker d\pi$). En un point $p\in Z(M,g)$, comme $\HH_p$ est
isomorphe à $T_{\pi(p)}M$ via $d\pi$, la distribution horizontale
hérite naturellement de la structure presque complexe induite par
 $p$. La somme de cette structure presque complexe et de la structure
 complexe sur les
fibres munit $Z(M,g)$ d'une structure presque complexe naturelle
notée $\JJ$.

Le résultat fondamental sur lequel repose la théorie des twisteurs
est que l'intégra\-bilité de $\JJ$ dépend de la courbure de $g$.
Plus exactement, notons  $R$ le tenseur de courbure défini pour
tout champ de vecteurs $X,Y\in TM$ par :

 $$R(X,Y)
=[\nabla_X,\nabla_Y]-\nabla_{[X,Y]}.$$

 Comme $M$ est orientée, l'opérateur
de Hodge induit la décomposition \hbox{$\bw^2TM=\bw^+\oplus\bw ^-$} de
sorte qu'en tant qu'endomorphisme de $\bw^2 TM$, on a la
décomposition du tenseur de courbure \cite{ST, Bes87} 
$$
R=\left[\begin{array}{cc}W^++\frac{s}{12}Id&B\\
B^\star&W^-+\frac{s}{12}Id\end{array}\right].
$$
L'opérateur $W=W^++W^-$ est l'opérateur de Weyl, $s$ est la
courbure scalaire de $g$, $B$ son tenseur de Ricci sans trace,
$B^\star$ son adjoint et $Id$ la matrice identité.

\quad\\
{\bf Théorème \cite{AHS78}.} Pour toute $4$-variété riemannienne
orientée $(M,g)$, la  structure presque complexe  $\JJ$
 est intégrable :

\begin{enumerate}
\item[a)] sur $Z^+(M,g)$ si et seulement si  $g$ est
anti-autoduale, c'est-à-dire $W^+=0$ ;

\item[b)] sur $Z^-(M,g)$ si et seulement si  $g$ est autoduale,
c'est-à-dire $W^-=0$.
\end{enumerate}

\quad\\
{\bf Remarque.} On identifiera toujours un endomorphisme
anti-symétrique $u$ de $TM$ au bi-vecteur $\phi(u)$ de $\bw^{2} TM
$ via
$$g(\phi(u),X\w Y)=g(uX,Y)
\quad\forall X,Y\in TM.
$$
Ainsi on peut voir $Z^\pm(M,g)$ comme un sous-ensemble de
$\bigwedge^\pm$.

\section{Espaces des twisteurs "généralisés" en dimension 4}
 \subsection{Construction.} Ici encore, on considère $(M,g)$ une 4-variété
riemannienne orientée connexe.
 La métrique riemannienne sur $TM$ se prolonge en
une métrique sur $\T M$ encore notée $g$. On notera $\mathcal
O_g(\T M)$ le fibré des endomorphismes de $\T M$ orthogonaux pour
la pseudo-métrique $<.,.>$ et pour la métrique $g$.

\quad\\
{\bf Définition.} Le fibré des twisteurs généralisés d'une
4-variété riemannienne orientée est le fibré $\pi : \ZZ (M,g)\lra
M$ des structures presque complexes généralisées sur $M$
compatibles avec la métrique $g$
$$
\mathcal Z(M,g)=\{u\in \mathcal O_g(\T M)/u^2=-Id \}.
$$
\quad\\
 Comme nous l'avons vu avec l'exemple 1, toute structure
complexe sur $M$ compatible avec $g$ définit une structure
complexe généralisée compatible avec $g$. On a donc une inclusion
naturelle $Z(M,g)\subset\mathcal Z(M,g)$.

Nous verrons dans la partie 3.4 que la structure presque complexe $\JJ$
sur $Z(M,g)$ se prolonge de façon naturelle en une structure
presque complexe généralisée sur tout  $\mathcal Z(M,g)$. Mais
avant cela nous voudrions comparer notre définition à deux autres
définitions que nous pouvons trouver dans la littérature.

\subsection{Différence avec les constructions de Davidov et Mushkarov.}
En 2006,  Davidov et Mushkarov \cite{DM1} ont  étudié le fibré
$\mathcal G\lra M$ de toutes les structures presque complexes
généralisées sur $M$ sans demander de compatibilité avec une
métrique riemannienne. Comme dans le cas classique, si on se donne
une connexion sans torsion sur $M$, alors $\mathcal G$ admet une
structure presque complexe généralisée. La condition
d'intégrabilité est alors très restrictive:

\quad\\
{\bf Théorème \cite{DM1}.} Soit $M$ une 2$n$-variété munie d'une
connexion sans torsion $\nabla$. La structure presque complexe
généralisée (naturelle) sur $\mathcal G$ est intégrable si et
seulement si
\begin{enumerate}
\item[(i)] $n=1$ ;

 \item[(ii)] $n\geq 2$ et $\nabla$ est une
connexion plate.
\end{enumerate}

\quad\\
 Contrairement à la construction d'Atiyah, Hitchin, Singer et
 contrairement à la nôtre, les fibres de $\mathcal G\lra M$ ne sont
  pas compactes. Par contre $\ZZ(M,g)$ est une sous-variété
   de $\mathcal G$, c'est même une rétraction de $\mathcal G$.
    D'autre part, comme nous le verrons dans la partie 3.4,
     la condition d'intégrabilité de $\mathbb J$ sur $\mathcal Z(M,g)$
      sera beaucoup moins restrictive que celle sur $\mathcal
      G$.\\

 L'année suivante, Davidov et Mushkarov  \cite{DM2} se sont intéressés à
une autre variété qui ressemble à $\mathcal Z(M,g)$. Pour
comprendre leur construction rappelons ce qu'est une structure
kählérienne généralisée.

\quad\\
{\bf Définition \cite{Gua2,Gua1}.} Une structure (presque)
kählérienne généralisée sur $\T M$ est une paire $(\mathcal
J_1,\mathcal J_2)$ de structures (presque) complexes généralisées
qui commutent et telles que la 2-forme définie par
$$
G(\mathcal X,\mathcal Y)=<\mathcal J_1 \mathcal X,\mathcal J_2
\mathcal Y>,\quad \forall \mathcal X,\mathcal Y\in \T M,
$$
soit définie positive.

\quad\\
{\bf Exemple.} Soit $(M,J,w,g)$ une structure kählérienne
classique, c'est-à-dire  une structure complexe $J$, une structure
symplectique $w$ et une métrique riemannienne $g$ telles qu'on ait
le diagramme commutatif
$$\xymatrix{
   TM\ar[rr]^{g}& & T^\star M \\
&TM\,.\ar[ul]^{J}\ar[ur]^{w}& }
$$
La métrique riemannienne $g$ sur $TM$ s'étend en une métrique sur
$\T M$. En identifiant $\T M$ et $\T^\star M$ grâce à la
pseudo-métrique $<.,.>$,  la métrique $g$ peut-être vue comme
un endomorphisme $G=\left(\begin{array}{cc} 0&g^{-1}\\
g&0\end{array}\right)$ de $\T M$. Comme $\mathcal J_J\mathcal
J_w=\mathcal J_w\mathcal J_J= -G$, la paire $(\mathcal
J_J,\mathcal J_w)$ est une structure kählérienne généralisée sur
$\T M$.

\quad\\
 Lorsque $M$ est une variété de dimension deux, Davidov et Mushkarov
 introduisent le fibré $\mathcal P\lra M$ des structures
presque kählériennes généralisées sur $M$. Une nouvelle fois, si
$M$ est munie d'une connexion sans torsion, alors $\mathcal P$
admet deux structures presque complexes généralisées naturelles
suivant qu'on privilégie $\mathcal J_1$ ou $\mathcal J_2$.

\quad\\
{\bf Théorème \cite{DM2}.} Soit $M$ une 2-variété munie d'une
connexion sans torsion. Les structures presque complexes
généralisées (naturelles) sur $\mathcal P$ sont intégrables si et
seulement si la connexion est plate.

\quad\\
La proposition suivante montre que la différence entre notre
approche et la leur, outre la dimension de $M$, est le fait
d'imposer une contrainte métrique.

\quad\\
{\bf Proposition 1.} Soit $G$ l'endomorphisme de $\T M$ associé à
la métrique $g$. Le fibré $\ZZ(M,g)$ est le fibré des structures
presque kählériennes généralisées sur $M$ compatibles avec $g$ :
$$
\ZZ(M,g)\simeq\Big\{(\mathcal J_1,\mathcal J_2)\textrm{ structures
presque kählériennes généralisées}/ \mathcal J_1\mathcal
J_2=\mathcal J_2\mathcal J_1=-G\Big\}.
$$

\quad\\
{\bf Preuve.} Comme $G^2=Id$, on note $C^\pm$ le sous-espace
propre de $G$ associé à la valeur propre $\pm1$. Si, au-dessus
d'un ouvert $\mathcal U$ de $M$, on se donne
$(\theta_1,\ldots,\theta_{4})$ une base orthonormée de $TM$ et si
on note $(\theta_1^\star,\ldots,\theta_{4}^\star)$ sa base duale,
alors
$$
C^\pm=Vect\Big(\theta_1\pm \theta_1^\star,\ldots, \theta_{4}\pm
\theta_{4}^\star\Big).
$$
Comme les éléments de $\mathcal Z(M,g)$ sont compatibles avec $g$,
ils stabilisent les espaces propres $C^\pm$. La matrice d'un
élément $u\in\mathcal Z(M,g)$ dans une base adaptée à la
décomposition $\T M=C^+\oplus C^-$ est donc de la forme:
$\left(\begin{array}{cc} u_1&0\\0&u_2\end{array}\right)$ avec
$u_1,u_2\in Z(\R^4,e)$, où $e$ désigne la métrique euclidienne
canonique sur $\R^4$.
De plus, tout élément $u=\left(\begin{array}{cc} u_1&0\\
0&u_2\end{array}\right)$ de $\ZZ(M,g)$ définit la structure
presque kählérienne  généralisée $(\J_1,\J_2)$ sur $M$ compatible
avec $g$ où
$$\J_1=u=\left(\begin{array}{cc} u_1&0\\
0&u_2\end{array}\right)\textrm{ et }\J_2=
\left(\begin{array}{cc} u_1&0\\
0&-u_2\end{array}\right).
$$
 Réciproquement, étant donnée une
structure presque kählérienne généralisée $(\mathcal J_1,\mathcal
J_2)$ sur $\T M$ compatible avec $g$, il existe \cite{Gua1} $u_1$
et $u_2$ deux structures presque complexes sur $TM$ compatibles
avec $g$ telles que dans la base $TM\oplus T^\star M$  on ait
$$
\mathcal J_1=\displaystyle\frac{1}{2}\left(\begin{array}{ll}u_1+u_2&u_1-u_2\\
u_1-u_2&u_1+u_2\end{array}\right),\quad
\mathcal J_2=\displaystyle\frac{1}{2}\left(\begin{array}{ll}u_1-u_2&u_1+u_2\\
u_1+u_2&u_1-u_1\end{array}\right).
$$
Dans une base $C^+\oplus C^-$ cela donne
$$
\mathcal J_1=\displaystyle\left(\begin{array}{ll}u_1&0\\
0&u_2\end{array}\right),\quad
\mathcal J_2=\displaystyle\left(\begin{array}{ll}u_1&0\\0&-u_2\\
\end{array}\right).\; \square
$$

\subsection{Composantes connexes.}  En conservant les notations de
la preuve précédente,  on voit que l'espace des twisteurs
$\ZZ(M,g)$ admet les 4 composantes connexes suivantes :
$$
\ZZ^{++}(M,g)\simeq\left\{\left(\begin{array}{cc} u_1&0\\
0&u_2\end{array}\right)/ (u_1,u_2)\in Z^+(\R^4,e)\times
Z^+(\R^4,e)\right\}
$$
$$
\ZZ^{--}(M,g)\simeq\left\{\left(\begin{array}{cc} u_1&0\\
0&u_2\end{array}\right)/ (u_1,u_2)\in Z^-(\R^4,e)\times
Z^-(\R^4,e)\right\}
$$$$
\ZZ^{+-}(M,g)\simeq\left\{\left(\begin{array}{cc} u_1&0\\
0&u_2\end{array}\right)/ (u_1,u_2)\in Z^+(\R^4,e)\times
Z^-(\R^4,e)\right\}
$$$$
\ZZ^{-+}(M,g)\simeq\left\{\left(\begin{array}{cc} u_1&0\\
0&u_2\end{array}\right)/ (u_1,u_2)\in Z^-(\R^4,e)\times
Z^+(\R^4,e)\right\}.
$$
 D'autre part on a vu que le
fibré des twisteurs (classiques) $Z(M,g)$ pouvait être considéré
comme un sous-fibré de
 $\ZZ(M,g)$, plus exactement
$$Z^+(M,g)\simeq
\left\{\left(\begin{array}{ll}u_1&0\\0&u_2\end{array}\right)\in
\ZZ^{++}(M,g)/u_1=u_2\right\},$$

$$Z^-(M,g)\simeq
\left\{\left(\begin{array}{ll}u_1&0\\0&u_2\end{array}\right)\in
\ZZ^{--}(M,g)/u_1=u_2\right\}.
$$
 L'orientation sur $TM$ induit une orientation sur $T^\star
M$ et donc sur $\T M$. Pour finir le parallélisme avec le cas
classique, il est naturel d'introduire le fibré des structures
presque complexes généralisées compatibles avec la métrique et
l'orientation, soit
$$\ZZ^+(M,g)=\{u\in \OO_g(\T M)/u^2=-Id\textrm{ et } u>>0\}$$
$$\ZZ^-(M,g)=\{u\in \OO_g(\T M)/u^2=-Id\textrm{ et } u<<0\}$$
et de voir le lien avec ce que nous avons déjà défini. On définit
la parité d'une structure complexe généralisée comme la parité de
son type. Dans le cas où la dimension de $M$ est un multiple de
$4$, on sait que les endomorphismes $u_1$ et $u_2$ qui définissent
une structure presque kählérienne généralisée doivent avoir la
même parité (paire, paire) ou (impaire, impaire). Dans le premier
cas $u_1$ et $u_2$ induisent la même orientation; dans le deuxième
ils induisent deux orientations opposées \cite{Gua2,Gua1}. Ce qui
donne

\quad\\
{\bf Proposition 2.}
$$\begin{array}{lll}\ZZ^+(M,g)&=&\ZZ^{++}(M,g)\sqcup \ZZ^{--}(M,g)\\
&=&\Big\{(\J_1,\J_2)/ \J_1\J_2=-G
 \textrm{ et le type de $\J_1$ et de $\J_2$ est pair}
\Big\}, \end{array}$$
$$\begin{array}{lll}
\ZZ^-(M,g)&=&\ZZ^{+-}(M,g)\sqcup \ZZ^{-+}(M,g)\\
&=&\Big\{(\J_1,\J_2)/ \J_1\J_2=-G
 \textrm{ et le type de $\J_1$ et de $\J_2$ est impair}
\Big\}. \end{array}$$

\subsection{Structure complexe généralisée sur $\ZZ(M,g)$.}
 Soit $\V=ker d\pi$ l'espace vertical tangent aux fibres de $\pi
: \ZZ(M,g)\lra M$. La connexion de Levi-Civita sur $M$ nous
fournit une distribution horizontale $\HH$ qui est en somme
directe avec $\V$ : $T \ZZ(M,g)=\HH\oplus\V$. On identifiera le
dual $\V^\star$ (resp $\HH^\star$) avec les formes sur
$\T\ZZ(M,g)$ nulle sur $\HH$ (resp. sur $\V$).

 Le
groupe structural des fibrés $\ZZ^{++}(M,g),\;
\ZZ^{-+}(M,g),\;\ZZ^{+-}(M,g),\;\ZZ^{--}(M,g)$
 est $SO(3)$ et leurs fibres
s'identifient à $\Sp^2\times\Sp^2$. Il existe donc une structure
complexe  sur les fibres de $\V$ et donc une structure complexe
généralisée sur les fibres de
 \hbox{$ \V\oplus \V^\star\lra \ZZ(M,g)$}. De plus en un point $p\in\ZZ(M,g)$ comme
$\HH_p\oplus\HH_p^\star$ est isomorphe à $\T_{\pi(p)}M$ via
$d\pi$, ce sous-espace hérite naturellement de la structure
presque complexe généralisée induite par $p$. La somme de cette
structure presque complexe généralisée et de celle sur $\V\oplus
\V^\star$ définit une structure presque complexe généralisée sur
$\ZZ(M,g)$, qui sera notée $\JJ$.

\quad\\
{\bf Remarque.} La structure presque complexe généralisée $\mathbb
J$ n'est pas la $B$-transformation d'une structure symplectique ni
d'une structure complexe sur $\ZZ(M,g)$. Il y a même "un phénomène
de saut" pour le type. Plus précisément la proposition 2 nous dit
 que $\mathbb J$ est :
\begin{enumerate}

\item[a)] de type trois sur $\ZZ(M,g)^{+-}$ et sur
$\ZZ^{-+}(M,g)$,

 \item[b)]  de type quatre sur $Z(M,g)\subset \ZZ^{++}(M,g)\cup
\ZZ^{--}(M,g)$,

\item[c)] de type deux sur le complémentaire de $Z(M,g)$ dans $
\ZZ^{++}(M,g)\cup \ZZ^{--}(M,g)$.

\end{enumerate}

\quad\\
 Le résultat principal de cet article est le suivant.

\quad\\
{\bf Théorème 1.} Pour toute $4$-variété riemannienne orientée
$(M,g)$, la structure presque complexe généralisée $\mathbb J$ sur
son espace de twisteurs généralisés $\ZZ(M,g)$ est intégrable :

\begin{enumerate}
\item[a)] sur $\ZZ^{++}(M,g)$ si et seulement si  $g$ est
anti-autoduale et Ricci plate;

\item[b)] sur $\ZZ^{--}(M,g)$ si et seulement si $g$ est autoduale
et Ricci plate;

\item[c)] sur $\ZZ^{+-}(M,g)$ si et seulement si  $g$ est  à
courbure sectionnelle constante;

\item[d)] sur $\ZZ^{-+}(M,g)$ si et seulement si  $g$ est  à
courbure sectionnelle constante.
\end{enumerate}

\quad\\
Les variétés à courbure sectionnelle constante sont les variétés
dont le tenseur de courbure $R:\bigwedge^2TM\lra\bigwedge^2TM$ est
 une homothétie. Par changement conforme on peut supposer que la
courbure sectionnelle est constante égale à -1, 0 ou 1. Les deux
propositions suivantes nous donne les 4-variétés compactes sur
lesquels on peut appliquer le théorème.

\quad\\
{\bf Proposition \cite{Lee}.} Si $(M,g)$ est une 4-variété
complète à courbure sectionnelle constante égale à -1, 0 ou 1,
alors $M$ est isométrique au quotient $\widetilde M/\Gamma$ où :
\begin{enumerate}
\item[a)] $\widetilde M$ est  la sphère $\Sp^4$, la plan  $\R^4$
ou l'espace hyperbolique $\mathbb H^4$ munis de leur métrique
usuelle,

\item[b)] $\Gamma$ est un sous-groupe discret du groupe des
isométries de $\widetilde M$ isomorphe au groupe fondamental de
$M$ et dont l'action sur $\widetilde M$ est propre, discontinue et
sans point fixe.
\end{enumerate}

\quad\\
{\bf Proposition \cite{Hit74}.} Soit $(M,g)$ une 4-variété
compacte orientée munie d'une métrique anti-autoduale et Ricci
plate. Alors
\begin{enumerate}
\item[a)] soit $(M,g)$ est plat,

\item[b)] soit le revêtement universel de $M$ est une surface
$K3$.
\end{enumerate}
Dans le deuxième cas il y a les surfaces K3, les surfaces
d'Enriques \cite{BHPVdV} et leurs quotients par une involution
anti-holomorphe.
 En particulier, si $M$ est une surface d'Enriques,  on a  une structure complexe
généralisée non triviale sur $\ZZ^{++}(M,g)$. Pour les surfaces K3
et les tores plats qui sont hyperkählériens, l'intégrabilité de
$\JJ$ sur $\ZZ^{++}(M,g)$ était un résultat déjà connu de
Bredthauer \cite{Bre}.

\section{Démonstration}
\subsection{Lemme technique.} Soit $(M,g)$ une
4-variété riemannienne, $\nabla$ la connexion de Levi-Civita sur
$M$ et $\ZZ(M,g)$ l'espace des twisteurs associés.
  Soit $\mathcal U$ un petit ouvert de $M$ sur lequel on a une
  trivialisation de $\pi:\ZZ(M,g)\lra M$ et $(m,u)$
   des coordonnées sur  $\pi^{-1}(\mathcal U)$.

   Nous noterons $\overrightarrow \X\in T\ZZ(M,g)$ la partie
vectorielle de $\X\in \T \ZZ(M,g)$, c'est-à-dire la projection sur
$T\ZZ(M,g)$ parallèlement à $T^\star \ZZ(M,g)$:
$$\begin{array}{lccc}
\overrightarrow{}:&\T \ZZ(M,g)&\lra& T\ZZ(M,g)\\
&\X=X+\xi&\lms&\overrightarrow{\X}:=X.
\end{array}
$$
La connexion de Levi-Civita s'étend sur $\T M$. On définit $\RR$
le tenseur de courbure de $g$ par
$$
\RR(X,Y)\mathcal Z=[\nabla_X,\nabla_Y]\mathcal Z
-\nabla_{[X,Y]}\mathcal Z, \quad \forall X,Y\in \mathcal C^\infty(
TM) \textrm{ et }\forall \mathcal Z\in \mathcal C^\infty(\T M).
$$
Pour alléger l'écriture, si $\X,\Y\in T \ZZ(M,g)$ sont des champs
de vecteurs sur $\ZZ(M,g)$ on écrira $\RR(\X,\Y)$ plutôt que
$\RR(\pi_\star\X,\pi_\star\Y)$.

\quad\\
 Pour démontrer le théorème 1 nous aurons besoin du lemme technique suivant.

\quad\\
{\bf Lemme  technique A.}  Soit $(M,g)$ une 4-variété
riemannienne. La structure presque complexe généralisée $\JJ$ sur
$\ZZ(M,g)$ est intégrable si et seulement si, pour tous champs
$\X,\Y\in \HH\oplus\HH^\star$ et pour tout point
$(m,u)\in\ZZ(M,g)$, on a
$$\left[u,\RR\Big(\overrightarrow{\X}\wedge
\overrightarrow{\Y}
  -\overrightarrow{u\X}\wedge
    \overrightarrow{u\Y}\Big)+
     u\RR\Big(\overrightarrow{u\X}\wedge\overrightarrow{\Y}
     +\overrightarrow{\X}\wedge
    \overrightarrow{u\Y}\Big)\right]=0.$$

\quad\\
 La suite de cette section est consacrée à la démonstration
de ce lemme. Soit $X+\xi$ une section de $\T
 M\lra M$. On notera $\wh X+\wh \xi\in\HH\oplus\HH^\star$ le champ
  relevé. Un tel champ est dit {\it basique}.

\quad\\
{\bf Proposition 3.} Soient  $A,B\in \V$ deux champs de vecteurs
verticaux sur $\ZZ(M,g)$. On se donne $X\in TM$ un champ de
vecteurs sur $M$ et $\xi \in T^\star M$ une forme sur $M$.
\begin{enumerate}

\item $[A,B]\in\V$ ,

\item $[\wh X,A]\in\V$ ,

\item $[\wh X+\wh \xi ,\JJ A]=\JJ[\wh X+\wh \xi,A]$ ,

\item $[\JJ(\wh X+\wh\xi),\JJ A]=\JJ[\JJ(\wh X+\wh \xi),A]$.

\end{enumerate}

\quad\\
{\bf Preuve.}  Le premier point provient du fait que la
distribution verticale est l'espace tangent aux fibres de $\pi :
\ZZ(M,g)\lra M$. Comme $\wh X$ est un champ relevé, le deuxième
point est immédiat. Et comme le transport parallèle suivant les
directions horizontales respecte l'orientation et la métrique sur
les fibres, il respecte la structure complexe sur l'espace tangent
vertical, on a $[\wh X,\JJ U]=\JJ[\wh X,U]$. De plus, par
définition du crochet de Courant, on vérifie que
$[\wh\xi,U]=0=[\wh\xi,\JJ U]$, ce qui termine la preuve du
troisième point. Le quatrième point s'obtient par "linéarité". En
effet soit $(\wh\X_1,\ldots,\wh\X_{8})\in \HH\oplus\HH^\star$ une
base de champs de vecteurs et de formes horizontales basiques.
Comme $\JJ$ stabilise $\HH\oplus\HH^\star$, on note $[\JJ_{ij}]$
la matrice de la restriction de $\JJ$ à $\HH\oplus\HH^\star$ dans
cette base. En utilisant le point 3, on a
$$
\begin{array}{ccl}
\JJ[\JJ\wh\X_j,A]&=&\JJ[\JJ_{ij}\wh\X_i,A]\\
&=&\JJ\big(\JJ_{ij}[\wh\X_i,A]-A\wh\X_j\big)\\
&=&\JJ_{ij}[\wh\X_i,\JJ A]-\JJ A\wh\X_j.
\end{array}
$$
et
$$
\begin{array}{ccl}
[\JJ\wh\X_j,\JJ A]&=&[\JJ_{ij}\wh\X_i,\JJ A]\\
&=&\JJ_{ij}[\wh\X_i,\JJ A]-\JJ A\wh\X_j.\; \square
\end{array}
$$

\quad\\
 {\bf Corollaire 1.} Le tenseur de Nijenhuis de $\JJ$ sur $\ZZ(M,g)$ vérifie
 $\NN(\X,A)=0$ pour tout
$\X\in\HH\oplus\HH^\star$ et pour tout $A\in\V$.

\quad\\
{\bf Preuve.} Par linéarité, on peut supposer que $\X$ est un
champ basique. Le corollaire 1 est alors une conséquence immédiate
des points 3 et 4 de la proposition 3. $\square$

\quad\\
{\bf Proposition 4.} Soient $X, Y\in TM$ deux champs de vecteurs
sur $M$. Au-dessus de l'ouvert $\mathcal U$,  la décomposition du
champ de vecteurs $[\wh X,\wh Y]$ en partie horizontale et
verticale au point $(m,u)$ est donnée par :
$$[\wh X,\wh Y] =\wh{[X,Y]}+[u,\RR(X,Y)].$$

\quad\\
{\bf Preuve.} Notons $G$ le groupe $\mathcal O_g(\R^4\oplus
\R^{4\star})=O(4,4)\cap O(8)$ ; c'est un groupe à quatre
composantes connexes. Notons également
 $\theta$ la 1-forme de connexion sur le $G$-fibré principal $\mathcal
O_g(\T M)$ associée à la connexion de Levi-Civita. Soient $X,Y\in
TM$ deux champs de vecteurs sur $M$ et $\widetilde X,\widetilde Y$
leurs relevés horizontaux dans $\mathcal O_g(\T M)$. Avec la convention de signe que nous avons choisie pour le tenseur de courbure, la
décomposition du champ de vecteurs $[\widetilde X,\widetilde Y]$
en parties horizontale et verticale est donnée par (cf. \cite{O},
\cite{Bes87} chap 9)
$$
[\widetilde X,\widetilde
Y]=\widetilde{[X,Y]}+(\theta\vert_{\V})^{-1}(\RR(X,Y)),$$
 où, par définition, $(\theta\vert_{\V})^{-1}(\RR(X,Y))$ est le
 champ de vecteurs vertical sur $\mathcal O_g(\T M)$ défini au point
 $p\in \mathcal  O_g(\T M)$ par
 $$
\frac{d}{dt}\vert_{t=0}\Big(p.\,exp(t\RR(X,Y))\Big)=p.\,\RR(X,Y).
$$
 Le groupe $G$ agit transitivement sur les fibres de
$\ZZ(M,g)\lra M$ et la variété $\ZZ(M,g)$ est le fibré associé de
fibres $\Sp^2\times\Sp^2$. Plus précisément, le groupe $G$ agit à
droite sur $\mathcal  O_g(\T M)\times(\Sp^2\times\Sp^2)$
$$
\begin{array}{ccc}
\mathcal  O_g(\T M)\times (\Sp^2\times\Sp^2)\times G&\lra&\mathcal  O_g(\T M)\times(\Sp^2\times\Sp^2)\\
(p,j,g)&\lms&(p.g,g^{-1}.j)=(p.g,gjg^{-1})\,,
\end{array}
$$
et $\ZZ(M,g)$ est le quotient de $\mathcal  O_g(\T
M)\times(\Sp^2\times\Sp^2)$ par cette action. On notera $\Pi$ la
projection :
$$
\begin{array}{cccc}
\Pi &: \mathcal O_g(\T M)\times(\Sp^2\times\Sp^2)&\lra&\ZZ(M,\DD)\\
&(m,p,j)&\lms&(m,u)=(m,p^{-1}jp)\,.
\end{array}
$$
Comme $d\Pi\Big(p.\RR(X,Y)\Big)=[u,\RR(X,Y)]$, on a bien
$$
[\wh X,\wh Y]=\wh{[X,Y]}+[u,\RR(X,Y)].\;\square
$$

\quad\\
{\bf Corollaire 2.} Soit $U^\sharp\in\V^\star$ une 1-forme
verticale,
 $\X\in \HH\oplus \HH^\star$ un champ de vecteurs et de formes
 horizontal. Au point $(m,u)\in\ZZ(M,g)$, le tenseur de Nijenhuis $\NN(U^\sharp,\X)$ est la
1-forme horizontale définie  pour tout champ de vecteurs
horizontal $\overrightarrow{\Y}\in\HH$ par
$$
 \NN(U^\sharp,\X)(\overrightarrow{\Y})=U^\sharp\left(\left[u,\RR\Big(\overrightarrow{\X}\wedge
  \overrightarrow{\Y}  -\overrightarrow{u\X}\wedge
    \overrightarrow{u\Y}\Big)+
     u\RR\Big(\overrightarrow{u\X}\wedge\overrightarrow{\Y}
     +\overrightarrow{\X}\wedge
    \overrightarrow{u\Y}\Big)\right]\right).
 $$

 \quad\\
 {\bf Preuve.} Par définition du crochet de Courant, on sait que
 $[U^\sharp,\X]=[U^\sharp,\overrightarrow{\X}]$ est une 1-forme.
 Soient $A\in\V$ et $\overrightarrow{\Y}\in \HH$ deux
 champs de vecteurs. On a,
au point $(m,u)\in\ZZ(M,g)$,
$$
\begin{array}{ccl}
[U^\sharp,\X](A+\overrightarrow{\Y})&=&dU^\sharp(\overrightarrow{\X},\overrightarrow{\Y}+A)\\
&=&\overrightarrow{\X}.U^\sharp(A)-U^\sharp([\overrightarrow{\X},\overrightarrow{\Y}+A])\\
&=&\overrightarrow{\X}.U^\sharp(A)-U^\sharp([\overrightarrow{\X},A])
-U^\sharp([u,\RR(\overrightarrow{\X},\overrightarrow{\Y})])\,.
\end{array}
$$
Le point 3 de la proposition 3 nous assure alors que $[\JJ
U^\sharp,\X](A)=\JJ[U^\sharp,\X](A).$ La partie verticale de la
1-forme $\NN(U^\sharp,\X)$ est donc nulle.
 Pour la partie horizontale, le calcul précédent
nous dit qu'au point $(m,u)$  on a
$$
\NN(U^\sharp,\X)(\overrightarrow{\Y})=U^\sharp\left(\left[u,\RR\Big(\overrightarrow{\X}\wedge
  \overrightarrow{\Y}  -\overrightarrow{u\X}\wedge
    \overrightarrow{u\Y}\Big)+
     u\RR\Big(\overrightarrow{u\X}\wedge\overrightarrow{\Y}
     +\overrightarrow{\X}\wedge
    \overrightarrow{u\Y}\Big)\right]\right). \;\square
 $$

\quad\\
{\bf Corollaire 3.} Au point $(m,u)\in\ZZ(M,g)$ et pour tous
champs de vecteurs et de formes horizontaux
$\X,\Y\in\HH\oplus\HH^\star$, on a
$$\NN(\X,\Y)=-\left[u,\RR\Big(\overrightarrow{\X}\wedge \overrightarrow{\Y}
  -\overrightarrow{u\X}\wedge
    \overrightarrow{u\Y}\Big)+
     u\RR\Big(\overrightarrow{u\X}\wedge\overrightarrow{\Y}
     +\overrightarrow{\X}\wedge
    \overrightarrow{u\Y}\Big)\right].
    $$

\quad\\
{\bf Preuve.}
  On note
$(\X_1,\ldots,\X_{8})$ une  base orthonormée de $\T M$ (pour la
métrique et la pseudo-métrique). La distribution $\mathcal
H\oplus\HH^\star$ est stable par $\JJ$. On notera $[\JJ_{ij}]$ sa
matrice dans la base relevée $(\wh{\X_1},\ldots,\wh{\X_{8}})$ :
$$\begin{array}{rll}
  \left[\mathbb J \widehat{\X_i},\mathbb
  J\widehat{\X_j}\right]&=&\overrightarrow{\JJ\wh{\X_i}}.(\JJ_{rj})\;\widehat{\X_r}
  -\overrightarrow{\JJ\wh{\X_j}}.(\JJ_{li})
\;\widehat{\X_l}+  \JJ_{li}\JJ_{rj}\left[\widehat{\X_l},\widehat{\X_r}\right]\\
&&-\JJ_{ri}d\JJ_{rj}+\JJ_{lj}d\JJ_{li}\\
  \,\left[\mathbb J\widehat{\X_i},\widehat{\X_j}\right]
  +\left[\widehat{\X_i},\mathbb
  J\widehat{\X_j}\right]  &=&-\overrightarrow{\wh{\X_j}}.(
  \JJ_{li})\;\widehat{\X_l}+\JJ_{li}\left[\widehat{\X_l},
  \widehat{\X_j}\right]+ \overrightarrow{\wh{\X_i}}.(
  \JJ_{rj})\;(\widehat{\X_r})
  +\JJ_{rj}\left[\widehat{\X_i},\widehat{\X_r}\right]\\
  &&+d\JJ_{ji}-d\JJ_{ij}.
\end{array}
$$
En utilisant la proposition 4, on en déduit qu'au point $(m,u)$,
la partie verticale de $\NN(\wh\X_i,\wh\X_j)$ vaut
$$
-\left[u,\RR\Big(\overrightarrow{\X_i}\wedge \overrightarrow{\X_j}
  -\overrightarrow{u\X_i}\wedge
    \overrightarrow{u\X_j}\Big)+
     u\RR\Big(\overrightarrow{u\X_i}\wedge\overrightarrow{\X_j}
     +\overrightarrow{\X_i}\wedge
    \overrightarrow{u\X_j}\Big)\right].
$$
Pour la partie horizontale, on se donne $s$ une section de
$\ZZ(M,g)\lra M$ telle que $s(m)=u$ et $\nabla_m s=0$. La partie
horizontale de $\NN(\wh\X_i,\wh\X_j)$ restreinte à la sous-variété
$s(M)$ est égale au relevé horizontal du tenseur de Nijenhuis
$\NN(\X_i,\X_j)$ de $M$ munie de la structure presque complexe
généralisée induite par $s$. Comme la connexion est sans torsion,
on en déduit qu'au point $m$ on a $\NN(\X_i,\X_j)=0$, donc la
partie horizontale de $\NN(\wh\X_i,\wh\X_j)$ est nulle au point
$(m,u)$ donc partout. $\square$

\quad\\
 Comme les fibres de
$\ZZ(M,g)\lra M$ sont complexes, il est clair que quels que soient
$A,B\in \V\oplus\V^\star$, on a $\NN(A,B)=0$. Le lemme technique A
est alors une conséquence directe des corollaires 1, 2 et 3.

\subsection{Démonstration du théorème 1.}
  Si on change l'orientation sur $M$, le fibré
$\ZZ^{--}(M,g)$ devient $\ZZ^{++}(M,g)$ tandis que $\ZZ^{-+}(M,g)$
devient $\ZZ^{+-}(M,g)$. Il suffit donc d'étudier l'intégrabilité
sur $\ZZ^{++}(M,g)$ et sur $\ZZ^{+-}(M,g)$.\\ La connexion de
Levi-Civita $\nabla$ stabilise $C^\pm$ ; donc, pour tous $X,Y\in
TM$, le tenseur de courbure $\RR(X,Y)$ aussi. Plus exactement, si $(\theta_
1,\ldots,\theta_4)$ est une base orthonormée de $TM$ et si $(\theta_1^\star,\ldots,\theta_4^\star)$ est sa base duale; alors la
matrice du tenseur de courbure $\RR(X,Y)$ dans la base $(\theta_1+\theta_1^\star,\ldots,\theta_4-\theta_4^\star)$, adaptée à la décomposition $\T M=C^+\oplus C^-$, s'écrit
$\left(\begin{array}{cc}R(X,Y)&0\\0&R(X,Y)\end{array}\right)$. De plus, dans cette base un élément $u\in \ZZ(M,g)$ est
de la forme
$\left(\begin{array}{cc}u_1&0\\0&u_2\end{array}\right)$. Suivant
que les vecteurs $\vec\X$ et $\vec\Y$ sont dans $C^+$ ou dans
$C^-$, la condition d'intégrabilité donnée par le lemme technique
A équivaut à l'annulation des six tenseurs définis pour tout triplet
$(m,u_1,u_2)\in\ZZ(M,g)$ et pour tous $X,Y\in TM$ par

\begin{enumerate}
\item[a)] $ G_1(X,Y,u_1,u_2)=\Big[u_1,\;\;R(X\wedge Y-u_1 X\wedge
u_1Y)+u_1 R(u_1X\wedge Y+X\wedge u_1Y)\Big]$,

\item[b)] $ G_2(X,Y,u_1,u_2)=\Big[u_2,\;\;R(X\wedge Y-u_1 X\wedge
u_1Y)+u_2 R(u_1X\wedge Y+X\wedge u_1Y)\Big]$,

\item[c)] $ G_3(X,Y,u_1,u_2)=\Big[u_1,\;\;R(X\wedge Y-u_2 X\wedge
u_2Y)+u_1 R(u_2X\wedge Y+X\wedge u_2Y)\Big]$,

\item[d)] $ G_4(X,Y,u_1,u_2)=\Big[u_2,\;\;R(X\wedge Y-u_2 X\wedge
u_2Y)+u_2 R(u_2X\wedge Y+X\wedge u_2Y)\Big]$,

\item[e)] $ G_5(X,Y,u_1,u_2)=\Big[u_1,\;\;R(X\wedge Y-u_1 X\wedge
u_2Y)+u_1 R(u_1X\wedge Y+X\wedge u_2Y)\Big]$,

\item[f)] $ G_6(X,Y,u_1,u_2)=\Big[u_2,\;\;R(X\wedge Y-u_1 X\wedge
u_2Y)+u_2 R(u_1X\wedge Y+X\wedge u_2Y)\Big]$.

\end{enumerate}
Sur un petit ouvert $\mathcal U$ de $M$, on se fixe
$(\theta_1,\theta_2,\theta_3,\theta_4)$ un champ de bases
orthonormées directes de $TM$. Cela définit une trivialisation
locale des fibrés $\bw^\pm\lra M$ via les sections
$$
\begin{array}{lll}
\left\{\begin{array}{lll}
I^+&=&\theta_1\w\theta_2+\theta_3\w\theta_4\\
J^+&=&\theta_1\w\theta_3-\theta_2\w\theta_4\\
K^+&=&\theta_1\w\theta_4+\theta_2\w\theta_3
\end{array}\right.&
\textrm{ et }& \left\{\begin{array}{llr}
I^-&=&\theta_1\w\theta_2-\theta_3\w\theta_4\;\\
J^-&=&\theta_1\w\theta_3+\theta_2\w\theta_4\;\\
K^-&=&-\theta_1\w\theta_4+\theta_2\w\theta_3.
\end{array}\right.
\end{array}
$$

\quad\\
 On rappelle qu'on identifie les éléments de $\bw^2TM$
aux endomorphismes anti-symétriques de $TM$. Un petit calcul
permet alors de vérifier que

\quad\\
{\bf Lemme 1.} Pour tout $(u_1,u_2)\in\bw^+\times\bw^-$ et pour
tous $X,Y\in TM$ on a

\begin{enumerate}
\item $X\w Y-u_1 X\w u_1 Y\in\bw^+$,

\item $u_1X\w Y+X\w u_1 Y\in\bw^+$ ,

\item $[u_1,u_2]=0$.
\end{enumerate}

\quad\\
\underline{\textsc{\'Etude de l'intégrabilité sur
$\ZZ^{++}(M,g)$}}

\quad\\
On a vu dans la partie 3.3 que, sur $\ZZ^{++}(M,g)$, les éléments
$u_1$ et $u_2$ variaient dans le même espace $Z^{+}(\R^4,e)$.
L'annulation des tenseurs $G_1,G_2,...,G_6$ équivaut donc à celle
des tenseurs $G_2$ et $G_6$. Sur la sous-variété $Z^+(M,g)\subset
\ZZ^{++}(M,g)$, l'annulation du tenseur $G_2$ est la contrainte
d'intégrabilité sur les espaces de twisteurs classiques, elle
entraîne donc que $g$ est anti-autoduale. L'annulation de $G_2$
impose également à la courbure scalaire d'être nulle car
$$
G_2(\theta_1,\theta_3,I^+,J^+)=\frac{1}{12}[J^+,s(J^++J^+K^+)]=-\frac{s}{6}K^+.
$$
D'autre part on a
$$
\begin{array}{lll}
G_6(\theta_1,\theta_1,I^+,J^+)&=&\displaystyle\frac{1}{2}\Big[J^+,-R(K^-)+J^+R(J^--I^-)\Big],\\
G_6(\theta_3,\theta_3,I^+,J^+)&=&\displaystyle\frac{1}{2}
\Big[J^+,\;\;\;R(K^-)+J^+R(J^-+I^-)\Big].
 \end{array}$$
  Si $G_6=0$, on doit donc avoir
$\Big[J^+,R(J^-)\Big]=0$ et de même
$\Big[I^+,R(J^-)\Big]=\Big[K^+,R(J^-)\Big]=0$. Soit $R(J^-)\in
\bigwedge^-$ et plus généralement $R(\bigwedge^-)\in \bigwedge^-$,
ce qui impose à la métrique $g$ d'être Einstein.

Réciproquement, si la  métrique $g$ est anti-autoduale et Ricci
plate, alors l'image de $R$ est dans $\bigwedge^-$, et  donc
$[\bw^+,R]=0$. D'où l'annulation des tenseurs $G_2$ et $G_6$.

\quad\\
\underline{\textsc{\'Etude de l'intégrabilité sur
$\ZZ^{+-}(M,g)$}}

\quad\\
L'annulation de $G_1$ (resp. de $G_4$) est la contrainte
d'intégrabilité sur l'espace des twisteurs $Z^+(M,g)$ (resp.
$Z^-(M,g)$). On a donc $G_1=G_4=0$ si et seulement si $g$ est
localement conformément plate. De plus sur $\ZZ^{+-}(M,g)$,
l'annulation de $G_2$ équivaut à demander à $g$ d'être d'Einstein.
En effet :
$$\begin{array}{rclc}
G_2(\theta_1,\theta_3,I^+,u_2)&=& \Big[u_2,\;\;
R(J^+)+u_2R(K^+)\Big]&\qquad (1)\\
G_2(\theta_1,\theta_3,K^+,u_2)&=&\Big[u_2,\;\;
R(J^+)-u_2R(I^+)\Big]& \qquad (2)\\
G_2(\theta_1,\theta_2,J^+,u_2)&=& \Big[u_2,\;\;
R(I^+)-u_2R(K^+)\Big]& \qquad (3)
\end{array}
$$
Si $G_2=0$, la combinaison $(1)-(2)+(3)$ donne  $
(Id+u_2)[u_2,R(I^+)]=0$, et donc $[u_2,R(I^+)]=0$ pour tout $
u_2\in \bigwedge^-$. On a donc $R(I^+)\in\bigwedge^+$ et plus
généralement $R(\bigwedge^+)\subset\bigwedge^+$ soit $B=0$.
Réciproquement, si $g$ est d'Einstein, le lemme 1 montre que
$G_2=0$. Si $\JJ$ est intégrable, nécessairement $g$ est à
courbure sectionnelle constante. Pour la réciproque, nous aurons
besoin du lemme suivant.

\quad\\
{\bf Lemme 2.} Soit $(M,g)$ est une 4-variété à courbure
sectionnelle constante, c'est-à-dire $R=\lambda Id$.  Sur notre
ouvert $\mathcal U$ et pour tous champs de vecteurs $X,Y\in TM$,
on a
$$
[I^+,R(X\w Y)]=\lambda\Big(- g(K^+X,Y)J^++g(J^+X,Y)K^+\Big).
$$

\quad\\
{\bf Preuve du lemme 2.}  Par symétrie et par bilinéarité, on peut
supposer que $X=\theta_1$ et $Y=\theta_i$
$$\begin{array}{lll}
R(\theta_1\w\theta_2)=\displaystyle\frac{\lambda(I^++I^-)}{2},&\textrm{
d'où }&
[I^+,R(\theta_1,\theta_2)]=0\, ,\\
R(\theta_1\w\theta_3)=\displaystyle\frac{\lambda(J^++J^-)}{2},&\textrm{
d'où }&
[I^+,R(\theta_1,\theta_3)]=\lambda K^+\, ,\\
R(\theta_1\w\theta_4)=\displaystyle\frac{\lambda(K^+-K^-)}{2},&\textrm{
d'où }& [I^+,R(\theta_1,\theta_4)]=-\lambda J^+\, ,
\end{array}
$$
ce qui donne directement le lemme 2. $\square$

\quad\\
Si $(M,g)$ est une variété à courbure sectionnelle constante, elle
est en particulier autoduale, anti-autoduale et d'Einstein, soit
$G_1=G_4=G_2=0$. Enfin, quitte à changer de base orthonormée, on
peut supposer que $u_1=I^+$.  D'après le lemme 2, on a alors
$$
\begin{array}{lll}
&&G_5(\theta_i,\theta_j,I^+,u_2)\\
&=&\quad\lambda\Big(-g\big(K^+\theta_i,\theta_j\big)+g\big(J^+\theta_i,u_2\theta_j\big)
+g\big(K^+\theta_i,\theta_j\big)-g\big(J^+\theta_i,u_2\theta_j\big)\Big)J^+\\
 &&+\;
\lambda\Big(g\big(J^+\theta_i,\theta_j\big)+g\big(K^+\theta_i,u_2\theta_j\big)
-g\big(J^+\theta_i,\theta_j\big)-g\big(K^+\theta_i,u_2\theta_j\big) \Big)K^+\\
&=&0.
\end{array}$$
Et en retournant l'orientation on a immédiatement $G_3=G_6=0$ : si
la métrique $g$ est à courbure sectionnelle constante, alors $\JJ$
est intégrable. $\square$

\section{Structure presque complexe sur $\ZZ(M,g)$}
\subsection{Critère d'intégrabilité.}
 Un des attraits de
la théorie des twisteurs est de faire un pont entre la géométrie
riemannienne et la géométrie complexe. Ainsi le résultat d'Atiyah,
Hitchin et Singer  traduit en terme d'intégrabilité d'une
structure complexe la propriété pour une métrique d'être
autoduale. Le théorème 1 exprime, lui, le fait d'être à courbure
sectionnelle constante ou autoduale et Ricci plate en terme
d'intégrabilité d'une structure complexe généralisée. On se
propose ici d'ajouter deux passerelles en caractérisant d'une part
les métriques d'Einstein et d'autre part les métriques autoduales
à courbure scalaire nulle.

 Pour cela, on  définit une
structure presque complexe sur $\ZZ(M,g)$. On sait déjà que la connexion de
Levi-Civita nous fournit une décomposition $T\ZZ(M,g)=\HH\oplus
\V$. En un point $(m,u)\in \pi^{-1}(\mathcal U)$ avec
$u=\left(\begin{array}{ll}u_1&0\\0&u_2\end{array}\right)$,  on
peut donc définir la structure presque complexe $\mathbb J_1$ sur $T
\ZZ(M,g)$ comme la somme directe de la structure complexe
naturelle sur $\V$  et  de l'action de $u_1$ sur $\HH$, après
identification avec  $TM$. Comme précédemment, l'intégrabilité de
cette structure complexe dépend de la courbure de $g$. Plus
exactement on a

\quad\\
{\bf Théorème 2.} Pour toute $4$-variété riemannienne orientée
$(M,g)$, la structure presque complexe $\mathbb J_1$ est
intégrable :

\begin{enumerate}
\item[a)] sur $\ZZ^{++}(M,g)$ si et seulement si  $g$ est
anti-autoduale à courbure scalaire nulle ;

\item[b)] sur $\ZZ^{--}(M,g)$ si et seulement si  $g$ est
autoduale à courbure scalaire nulle ;

\item[c)] sur $\ZZ^{+-}(M,g)$ si et seulement si    $g$ est
anti-autoduale et d'Einstein ;

\item[d)] sur $\ZZ^{-+}(M,g)$ si et seulement si  $g$ est
autoduale et d'Einstein.

\end{enumerate}

\quad\\
 {\bf Preuve.} Pour établir l'intégrabilité de $\mathbb
J_1$ on utilise le lemme technique suivant dont la démonstration
est similaire à celle du précédent.

\quad\\
{\bf Lemme technique B.} Au point $(m,u)\in\pi^{-1}(\mathcal U)$
avec $u=\left(\begin{array}{ll}u_1&0\\0&u_2\end{array}\right)$, le
tenseur de Nijenhuis de $\mathbb J_1$ est nul si et seulement si,
pour tout $u\in\ZZ(M,g)$ et pour tous  $X,Y\in TM$
$$\left[u,\RR\Big(\overrightarrow{X}\wedge
\overrightarrow{Y}
  -\overrightarrow{u_1X}\wedge
    \overrightarrow{u_1Y}\Big)+
     u\RR\Big(\overrightarrow{u_1X}\wedge\overrightarrow{Y}
     +\overrightarrow{X}\wedge
    \overrightarrow{u_1Y}\Big)\right]=0.
$$
\quad\\
 {\bf Preuve du théorème 2.} La
condition d'intégrabilité qui apparaît dans le lemme technique B
équivaut uniquement à l'annulation des deux tenseurs $G_1$ et
$G_2$ définis précédemment. Sur $\ZZ^{++}(M,g)$ on a vu que cela
  entraînait  $W^+=s=0$. Réciproquement, si $W^+=s=0$ alors $R(\bw^+)\subset \bw^-$
  et le lemme 1 nous donne automatiquement l'annulation des tenseurs
    $G_1$ et $G_2$.
     De m\^eme sur
$\ZZ^{+-}(M,g)$, on a vu que l'annulation de $G_1$ et de $G_2$
équivaut à $g$  anti-autoduale et d'Einstein. Enfin,
l'intégrabilité sur $\ZZ^{--}(M,g)$ et sur $\ZZ^{-+}(M,g)$ se
déduit de celle de $\ZZ^{++}(M,g)$ et de $\ZZ^{+-}(M,g)$ en
renversant l'orientation. $\square$

\subsection{Critère de semi-intégrabilité.} Les éléments de
$\ZZ(M,g)$ stabilisent $C^+$ et $C^-$. On introduit donc
$Z(C^-,g)$ l'ensemble des structures presque complexes sur $C^-$
compatibles avec $g\vert_{C^-}$, soit
$$
Z(C^-,g)=\{u\in End(C^-)/u^2=-Id \textrm{ et $u$ orthogonal pour
$g\vert_{C^-}$}\}.
$$
 Par restriction, on a la
projection naturelle de $\ZZ(M,g)$ sur $Z(C^-,g)$
$$
\begin{array}{cccc}
pr_- :& \ZZ(M,g)&\lra& Z(C^-,g)\\
&u&\lms&u\vert_{C^-}=u_2\,.
\end{array}
$$
{\bf Définition.} On dira qu'une structure presque complexe sur
$\ZZ(M,g)$ est {\it semi-intégrable} si la projection  sur
$Z(C^-,g)$ de son tenseur de Nijenhuis est nulle.

\quad\\
Techniquement, cela revient à demander l'annulation du tenseur
$G_2$. Sur $\ZZ^{++}(M,g)$, la semi-intégrabilité de $\JJ_1$
entraîne l'intégrabilité. Ce n'est pas le cas sur $\ZZ^{+-}(M,g)$,
où l'annulation de $G_2$ équivaut à demander à la métrique d'être
d'Einstein.

\quad\\
{\bf Proposition 5.} La structure presque complexe $\mathbb J_1$
est semi-intégrable sur $\ZZ^{+-}(M,g)$ (ou sur $\ZZ^{-+}(M,g)$)
si et seulement si la métrique $g$ est d'Einstein.

\quad\\
Cette proposition donne une interprétation complexe de la
propriété d'être Einstein.

\section{Dimension quelconque}

 La théorie des twisteurs
classiques comme nous l'avons présentée en dimension 4 a été
étendue en toute dimension \cite{Ber,O'br, DV,Sal, Leb, Slu} et en
particulier aux variétés quaternioniques Kähler.   Considérons
donc $(M,g)$ une $4n$-variété riemannienne avec $n>1$. Une
structure {\it presque hypercomplexe} sur $(M,g)$ est un triplet
$(I,J,K)$ de structures presque complexes compatibles avec $g$ et
telles que $IJ=-JI=K$. Lorsque $I,J,K$ sont intégrables, on parle
d'une structure {\it hypercomplexe}. Une structure {\it presque
quaternionique} sur $(M,g)$ est un sous-fibré de rang trois
\hbox{$D\subset End(TM)$} localement engendré par une structure
presque hypercomplexe. Une telle structure est dite {\it
quaternionique Kähler} si le fibré $D$ est préservé par la
connexion de Levi-Civita \cite{Bes87}.

Précisément, dans le cas où $(M,g,D)$ est une variété presque
quaternionique, on peut définir le fibré des twisteurs classiques
$Z(M,D)\lra M$ comme le fibré des structures presque complexes sur
$M$ appartenant à $D$. C'est un fibré de fibres $\Sp^2$.  Là
encore, on peut munir $Z(M,D)$ d'une structure presque
 complexe naturelle. Dans ce cas l'intégrabilité est automatique.

\quad\\
{\bf Théorème \cite{Bes87, Sal}.} Pour toute $4n$-variété
quaternionique Kähler $(M,g,D)$ avec $n>1$, la structure presque
complexe naturelle sur $Z(M,D)$ est toujours intégrable.

\quad\\
 Comme précédemment,  on se propose de donner une version
généralisée de ce théorème. Pour cela on considère $(M,g,D_1,D_2)$
une variété riemannienne munie de deux structures quaternioniques
Kähler. On note encore $C^\pm$ les sous-espaces propres de la
métrique $g$ vue comme endomorphisme de $\T M$. La projection de
$C^+$ (resp. $C^-$) sur $TM$ est un isomorphisme qui permet de
relever la distribution $D_1$ (resp. $D_2$) en une distribution
notée $D^+\subset End(C^+)$ (resp. $D^-\subset End(C^-)$).

\quad\\
{\bf Définition.} On définit $\ZZ(M,D_1,D_2)$ le fibré des
structures presque complexes compatibles avec $g$ dont la
restriction à $C^\pm$ appartient à $D^\pm$. Dans la base
$C^+\oplus C^-$ cela donne
$$
\ZZ(M,D_1,D_2)\simeq\left\{\left(\begin{array}{ll}u_1&0\\0&u_2\end{array}\right)/u_1\in
D^+ \textrm{ et } u_2\in D^-\right\}.
$$
 C'est un fibré de fibres $\Sp^2\times\Sp^2$.

\quad\\
 Le fibré des twisteurs classiques
 $Z(M,D)$
se réalise naturellement comme la sous-variété de $\ZZ(M,D,D)$
d'équation $u_1=u_2$. D'autre part, la même construction qu'en
dimension quatre permet de munir $\ZZ(M,D_1,D_2)$ d'une structure
presque
 complexe généralisée naturelle $\JJ$,  qui n'est pas la B-transformation
  d'une structure complexe,  ni d'une structure symplectique.

\quad\\
{\bf Remarque.} En dimension 4, $\bigwedge^+$ et $\bigwedge^-$
sont les seules structures presque quaternioniques sur $(M,g)$.
Par convention, on dira que la 4-variété $(M,g,\bigwedge^+)$ est
quaternionique Kähler si elle est anti-autoduale et d'Einstein.
Avec les notations de la section précédente on a
$$
\begin{array}{l}
\ZZ^{++}(M,g)=\ZZ(M,\bigwedge^+,\bigwedge^+)\,,\\
\ZZ^{+-}(M,g)=\ZZ(M,\bigwedge^+,\bigwedge^-)\,,\\
\ZZ^{-+}(M,g)=\ZZ(M,\bigwedge^-,\bigwedge^+)\,,\\
\ZZ^{--}(M,g)=\ZZ(M,\bigwedge^-,\bigwedge^-)\,.
\end{array}
$$

\quad\\
Avec cette convention, le théorème 1 se généralise à toute
dimension sous la forme suivante

\quad\\
{\bf Théorème 3.} Soit $n\geq 1$ et $(M,g,D_1,D_2)$ une
$4n$-variété munie de deux structures quaternioniques Kähler. La
structure presque complexe généralisée $\mathbb J$ sur
$\ZZ(M,D_1,D_2)$
  est intégrable si et seulement si l'une des deux conditions suivantes est
 vérifiée:
 \begin{enumerate}
 \item[a)] la courbure scalaire de $g$ est nulle,

 \item[b)] les éléments de $D_1$ et de $D_2$ commutent.
 \end{enumerate}

\quad\\
 Les variétés hyperkählériennes généralisées admettent naturellement
  deux structures quaternioniques Kähler sous-jacentes.
   Notre théorème s'applique donc en particulier à ces variétés, et
    on retrouve alors le résultat de Bredthauer \cite{Bre}.

\quad\\
{\bf Preuve.} La cas $n=1$ ayant déjà été traité, on suppose que
 $n>1$. Alors,  toutes les constructions précédentes marchent encore
 et, si on conserve  les mêmes notations, le théorème 3 est une
conséquence des trois lemmes suivants.

\quad\\
{\bf Lemme technique C.} Soit $(M,g,D_1,D_2)$ une $4n$-variété
munie de deux structures quaternioniques Kähler. La structure
presque complexe généralisée $\JJ$ sur $\ZZ(M,D_1,D_2)$ est
intégrable si et seulement si, pour tous champs $\X,\Y\in
\HH\oplus\HH^\star$ et pour tout point $u\in \ZZ(M,D_1,D_2)$, on a

$$\left[u,\RR\Big(\overrightarrow{\X}\wedge
\overrightarrow{\Y}
  -\overrightarrow{u\X}\wedge
    \overrightarrow{u\Y}\Big)+
     u\RR\Big(\overrightarrow{u\X}\wedge\overrightarrow{\Y}
     +\overrightarrow{\X}\wedge
    \overrightarrow{u\Y}\Big)\right]=0.$$

\quad\\
La démonstration de ce lemme est identique à celle du premier
lemme technique.

\quad\\
{\bf Lemma 3~\cite{Bes87}.}  Soit $n>1$, $(M,g,D)$ une
$4n$-variété quaternionique Kähler et $r$ le tenseur de Ricci de
la métrique $g$. On note $(I,J,K)$ une structure presque
hypercomplexe qui engendre $D$ au-dessus d'un ouvert de $M$.
Au-dessus de cet ouvert et, pour tous champs de vecteurs $X,Y\in
TM$, on a
  \[
  \begin{array}{rll}
    \,[I,R(X,Y)]&=&-\gamma(X,Y)J+\beta(X,Y) K\\
    \,[J,R(X,Y)]&=&-\alpha(X,Y)K+\gamma(X,Y)I\\
    \,[K,R(X,Y)]&=&-\beta(X,Y)I+\alpha(X,Y)J\\
  \end{array}
  \textrm{ avec } \left\{
    \begin{array}{l}
      \alpha(X,Y)=\frac{2}{n+2}r(IX,X)\\
      \beta(X,Y)=\frac{2}{n+2}r(JX,X)\\
      \gamma(X,Y)=\frac{2}{n+2}r(KX,X)\,.
    \end{array}
  \right.
  \]

\quad\\
{\bf Lemma 4~\cite{Bes87}.} Les variétés quaternioniques Kähler
sont d'Einstein.

\quad\\
Le lemme technique C nous assure que l'intégrabilité de $\JJ$ sur
$\ZZ(M,D_1,D_2)$ équivaut à l'annulation des six tenseurs
$G_1,\ldots, G_6$ définis dans la section 4.2. Si la courbure
scalaire de $g$ est nulle alors la métrique est Ricci plate et le
lemme 3 nous assure que ces six tenseurs sont nuls : $\JJ$ est
intégrable.

 On suppose donc que la courbure scalaire de $g$
est non nulle. On se fixe $\mathcal U$ un ouvert de $M$ et
$(I_1,J_1,K_1)$ (resp. $(I_2,J_2,K_2)$)
  une structure presque hypercomplexe
 sur $\mathcal U$ qui engendre $D_1$ (resp. $D_2$). Le lemme 3 nous dit
 que

 $$
 \begin{array}{ll}
&G_2(X,Y,I_1,I_2)=0\\
&\\
 \Longleftrightarrow& \Big( g(K_2X,Y)-g(K_2I_1X,I_1Y)+
 g(J_2I_1X,Y)+g(J_2X,I_1Y)\Big)
 J_2\\
 &+\Big(-g(J_2X,Y)+g(J_2I_1X,I_1Y)+
 g(K_2I_1X,Y)+g(K_2X,I_1Y)\Big)
 K_2=0\\
 &\\
 \Longleftrightarrow&\left\{\begin{array}{ll}
 K_2+I_1K_2I_1+[J_2,I_1]=0\\
 -J_2-I_1J_2I_1+[K_2,I_1]=0
  \end{array}\right.\\
 &\\
 \Longleftrightarrow&
 I_1[J_2,I_1]=[K_2,I_1].
 \end{array}
 $$
Et par symétrie $\left\{\begin{array}{l}\;I_1[I_2,I_1]=[J_2,I_1]\\
\;I_1[K_2,I_1]=[I_2,I_1]\end{array}\right.$, d'où l'on en déduit
que $[I_1,I_2]=0$ et plus généralement que les éléments de $D_1$
commutent avec ceux de $D_2$. La réciproque se vérifie facilement
en utilisant le lemme 3. $\square$

\quad\\
Guillaume DESCHAMPS,\\
\quad\\
 Université de Brest, UMR 6205,\\ Laboratoire de
Mathématiques de
Bretagne Atlantique\\
6 avenue Victor le Gorgeu\\
CS 93837,\\
29238 Brest cedex 3\\
(France)


\begin{thebibliography}{10}

\bibitem{AHS78}
M.F. Atiyah, N.J. Hitchin, I.M. Singer.
\newblock {\em Self-duality in four-dimensional {R}iemannian
geometry.}
\newblock {Proc. Roy. Soc. London Ser. A}, {\bf 362}, 425--461 (1978).


\bibitem{BHPVdV}
W.  Barth, K. Hulek, C.  Peters,
              A. Van de Ven.
\newblock {\em Compact complex surfaces.}
\newblock {Ergebnisse der Mathematik und ihrer Grenzgebiete. 3. Folge. A
              Series of Modern Surveys in Mathematics},
\newblock Springer-Verlag, Berlin (2004).





\bibitem{Ber}
L. Bérard-Bergery, T.  Ochiai.
\newblock {\em  On some generalizations of the construction of twistor spaces.}
\newblock Global Riemannian geometry (Durham, 1983),  Ellis Horwood Ser.
Math. Appl., Horwood, Chichester, 52-59, (1984).



\bibitem{Bes87}
A.L. Besse.
\newblock {\em Einstein manifolds.}  Ergebnisse der Mathematik
  und ihrer Grenzgebiete (3), {\bf 10}.
\newblock Springer-Verlag, Berlin (1987).



\bibitem{Bre}
A. Bredthauer.
\newblock {\em Generalized Hyperkähler geometry and supersymmetry.}
\newblock   Nucl. Phys. B, {\bf 773}, 172-183 (2007).


\bibitem{GC}
G. Cavalcanti, M. Gualtieri.
\newblock {\em A surgery for generalized complex structures on 4-manifolds.}
\newblock J. Differential Geom. {\bf 76}, 35-43, (2007).



\bibitem{Cou}
T. J. Courant.
\newblock {\em Dirac manifolds.}
\newblock   Trans. Amer. Math. Soc. {\bf 319}, 631-661 (1990).




\bibitem{DM1}
J. Davidov, O. Mushkarov.
\newblock {\em  Twistor spaces of generalized complex structures.}
\newblock J. Geom. Phys. {\bf 56}, 1623-1636 (2006).



\bibitem{DM2}
J. Davidov, O. Mushkarov.
\newblock {\em  Twistorial construction of generalized Kähler manifolds.}
\newblock  J. Geom. Phys. {\bf 57}, 889-901 (2007).


\bibitem{DV}
M. Dubois-Violette.
\newblock {\em  Structures complexes au-dessus des variétés,
applications.}
\newblock Mathematics and Physics,  Progr. Math. {\bf 37}, 1-42 (1983).





\bibitem{Gua2}
M. Gualtieri.
\newblock {\em Generalized complex geometry.}
\newblock Ph.D. thesis, St John's college, University of Oxford,
   arXiv: math/0401221, 107 pages (2003).


\bibitem{Gua1}
 M. Gualtieri.
\newblock {\em Generalized Kähler geometry.}
\newblock Commun. Math. Phys. {\bf 331}, 297-331 (2014).


\bibitem{Hit}
N.J. Hitchin.
\newblock {\em Generalized Calabi-Yau manifolds.}
\newblock  Q. J. Math. {\bf 54}, 281-308 (2003).


\bibitem{Hit74}
N.J. Hitchin.
\newblock {\em Compact four-dimensional Einstein manifolds.}
\newblock J. Differential Geom., {\bf 9}, 435-441
(1974).






\bibitem{Leb}
C. LeBrun.
\newblock {\em Quaternionic-Kähler manifolds and conformal
geometry.}
\newblock Math. Ann. {\bf 284}, 353-376 (1989).





\bibitem{Lee}
J.M. Lee.
\newblock {\em Riemannian manifolds. An introduction to curvature.}
Graduate Texts  Math., {\bf 176}.
\newblock Springer-Verlag, New York (1997).





\bibitem{O}
B. O'Neil
\newblock {\em The fundamental equations of a submersion.}
\newblock  Michigan Math. J. {\bf 13}, 459-469 (1966).







\bibitem{O'br}
N.R. O'Brian, J.H. Rawnsley.
\newblock {\em   Twistor spaces.}
\newblock  Ann. Global Anal.
Geom. {\bf 3}, 29-58  (1985).


\bibitem{Pen}
R. Penrose,  R.S. Ward.
\newblock {\em  Twistors for flat and curved space-time.}
\newblock General Relativity and Grav. {\bf 2}, 283-328, Plenum,
New York-London (1980).






\bibitem{Sal}
S. Salamon.
\newblock {\em Quaternionic Kähler manifolds.}
\newblock  Invent. Math. {\bf 67}, 143-171 (1982).



\bibitem{ST}
I.M. Singer, J.A. Thorpe.
\newblock {\em The curvature of {$4$}-dimensional Einstein spaces.}
\newblock Global Analysis, paper in honor of K. Kodaira, Princeton University Press,
 355-365 (1969).

\bibitem{Slu}
M. Slupinski.
\newblock {\em Espaces de twisteurs kählériens en dimension 4k, k>1.}
\newblock J. London Math. Soc. {\bf 33}, 535-542 (1986).







\end{thebibliography}
\end{document}